\documentclass[12pt]{article}
\usepackage[cp1251]{inputenc}
\usepackage[T2A]{fontenc}
\usepackage{epsfig}

\usepackage{a4wide}
\usepackage{amsmath}
\usepackage{amssymb}
\usepackage{amssymb,bbm,verbatim}
\usepackage{graphicx}
\usepackage{xcolor}

\newcommand{\ba}{\begin{eqnarray}}
\newcommand{\ea}{\end{eqnarray}}
\newcommand{\ban}{\begin{eqnarray*}}
\newcommand{\ean}{\end{eqnarray*}}

\newtheorem{theo}{Theorem}
\newtheorem{lem}[theo]{Lemma}
\newtheorem {prop}[theo]{Proposition}

\newtheorem {defi} {Definition}

\newcommand{\supp}{\mathop{\rm supp}\nolimits}

\newcommand{\abcd}[2]{\hbox to\textwidth{#1\dotfill #2}}

\renewcommand{\lll}{\left(}
\newcommand{\rrr}{\right)}

\newcommand{\md}{\mathop{\rm mod}\nolimits}

\begin{document}
\setcounter{page}{1}

\newpage
\begin{center}

{\Large {\bf  Harmonic Analysis on the Space of  $M$-positive Vectors}}
		
		\bigskip
		
		 {\large{\bf  Yu. Farkov$^1$ and \ M. Skopina$^{2.3}$}}
		
		\medskip  
		$^{1}$\ Russian Presidential Academy of National Economy and Public Administration,
		Moscow, 119571 Russia; farkov-ya@ranepa.ru
		\\
		$^{2}$\ St. Petersburg State University,  St. Petersburg, 199034 Russia
		\\
		$^{3}$\ Regional Mathematical Center of Southern Federal University,
		\\ 
		Rostov-on-Don, 344006 Russia;  skopinama@gmail.com

\end{center}

\bigskip

{\bf Abstract.} Given a dilation matrix $M$, a so-called space of $M$-positive vectors in the Euclidean space is introduced and studied. An algebraic structure of this space is similar to the positive half-line equipped with the termwise addition modulo 2, which is used in the Walsh analysis. The role of harmonics  is played by some analogues of the classical Walsh functions. 
The concept of Fourier transform is introduced, and the
Poisson summation formula, Plancherel theorem, Vilenkin-Chrestenson formulas and so on
are proved. A kind of analogue of the Schwartz class is studied. 
This class consists of functions such that both the function itself and its Fourier transform have compact support.

{\bf Keywords:} space of  $M$-positive vectors,  fractal sets, characters, Walsh functions, Fourier transform, step-functions.

\bigskip

{\large{\bf \S~1. Introduction}
	
	\medskip
	
The space of $M$-positive vectors $X=X(M, D)\subset {\mathbb R}^d$,
associated with a dilation matrix $M$ and a set of its digits $D$, is studied. This space $X$ is an analog of the positive half-line ${\mathbb R_+}$ equipped with the termwise addition modulo 2,  which is used in the Walsh harmonic analysis (see, e.g., \cite{SWS}-\cite{FS23}), and a  self-similar set $U=U(M,D)\subset X$ plays the role of the segment $[0,1]$. However, unlike the half-line,  the set $U$ as well as the whole space $X$, has a complex fractal structure for many matrices $M$.	
On the other hand, the space $X$ is  close to the Vilenkin groups in many aspects. In particular, the functions $\chi(\cdot,\omega)$ playing the role of characters in the Vilenkin groups, are very important. The Walsh functions for the space  $L^2(U)$  are defined on the basis of these functions,  and an analog of the Vilenkin-Chrestenson transform  is  presented.  
Also, by analogy with the Vilenkin groups, on the basis of the functions $\chi(\cdot,\omega)$, the concept of  Fourier transform is introduced. The basic properties of the Fourier transform are studied, in particular, the Plancherel theorem, the inversion formula and the  Poisson summation formula are proved. 
Thus, the harmonic analysis on the space $X$  looks quite usual in general. 
However one can see something  unusual for some classes of functions. In many problems, it is of great interest to have the maximum possible decay for both the function itself and its Fourier transform. In contrast to the classical harmonic analysis on ${\mathbb R}^d$, where  the function  and its Fourier transform cannot be compactly supported at the same time, such a class of functions exists in the space of $M$-positive vectors. This class consists of  step functions.
It  is proved  that the Fourier transform of a compactly supported  step function is also a compactly supported  step function. So, again we have  an analogy with the Vilenkin groups (see, e.g.,	\cite{GES, KP20}).

\bigskip

{\large{\bf \S~2. Space of $M$-positive vectors}

\medskip
As usual,   ${\mathbb R}^d$ is the Euclidean space with the inner product $\langle \cdot\ ,\cdot \rangle$, $\mu$ denotes the Lebesgue measure  on ${\mathbb R}^d$,  and  ${\mathbb Z}^d$ is the integer lattice of   ${\mathbb R}^d$. The characteristic function of a set $E\subset {\mathbb R}^d $ is denoted by $\mathbf{1}_E$.

Let $M\,$  be  an $d\times d$  integer matrix such that  $m:=|\det M|\geq 2$. 
Vectors $l,k\in\mathbb{Z}^d$ are said to be {\it congruent modulo} $M$  (we write $l \equiv k\,({\rm mod}\, M)$) if $l-k = Ms$ for some $s\in\mathbb{Z}^d$. 
The lattice ${\mathbb Z}^d$ is partitioned into cosets with respect  to this congruence. A set $D=D(M)$   is said to be a   {\it set of digits} for~$M$, if 
$
D=\{s_0, \,s_1, \dots, s_{m-1}\},
$
where   $s_0=\mathbf{0}\,$ is the zero vector,  $s_i\in\mathbb{Z}^d$ and  $s_i\not\equiv s_j \,({\rm mod}\, M)$ whenever $i\neq j$.  Let $M^*$ be the matrix transposed to $M$, and let the vectors
 $s^*_0, \,s^*_1, \dots, s^*_{m-1}$, where $s^*_0=\mathbf{0}$, form a set $D^*$ of digits for~ $M^*$, i.e.. $D^* = D(M^*)=\{s^*_0, \,s^*_1, \dots, s^*_{m-1}\}$. 
{\sloppy

}

\medskip
The following statement is proved in \cite[\S~2.2]{NPS}.

\begin{prop}\label{p1.2}
Let $M\,$ be an integer $d\times d$ matrix, $m=|\det M|\geq 2$, $D\,$ and $D^*\,$ be  sets of digits for  $M$ and  $M^*$, respectively. Then the following equality   
\begin{equation}
\sum_{s^*\in D^*}e^{ 2\pi i\langle M^{-1}l,\,s^* \rangle}=\left\{
\begin{array}{l}
\  m,\quad\mbox{if} \quad l = \mathbf{0}\,(\md M),\\
\  \, 0,\, \quad \mbox{if} \quad l \neq \mathbf{0}\,(\md M)
\end{array}
\right.
\label{1.1}
\end{equation}
holds for every $l\in\mathbb{Z}^d$.

\end{prop}

\medskip
In what follows, we will consider a class  $\mathfrak{M}_d$ of dilation matrices of
 order~$d$.
An integer  matrix $M$ is called {\it a dilation matrix}, if all its eigenvalues are bigger than~1 in absolute value.
For any matrix $M\in\mathfrak{M}_d$, there holds 
\ba
\label{my000}
\lim_{n\to +\infty}\|M^{-n}\|=0
\ea
and, moreover, the following relation 
$$
\sum_{n=1}^{\infty}\|M^{-n}\|^{\delta} < \infty,
$$
holds true  for any  $\delta>0$ (see, e.g., \cite[\S~2.2]{NPS}).

\medskip
Let  $D\,$ be a set of digits for a matrix  $M\in \mathfrak{M}_d$, and let $X^+=X^+(M,D)$ denote the set of vectors 
$x\in\mathbb{R}^d$ represented in the form
\begin{equation}
	x =  \sum_{j=-\infty}^{\infty}M^{-j}x_j, \quad x_{j} \in D,
	\label{1.81}
\end{equation}
where only a finite number of $x_j$ with a negative $j$ may be nonzero.
Let  $X^0=X^0(M,D)$ denote the set of vectors  $x\in X^+$ for which there exists a {\it finite representation},  i.e.,  expansion~\eqref{1.81} with  a finite number of nonzero~$x_j$.

\begin{prop}\label{p1.20}
For every $x\in X^0$ there exists a unique finite representation~\eqref{1.81}. 

\end{prop}

\medskip
{\bf Proof.} Assume that there exists a vector  $x\in X^0$ such that
$$
  x =  \sum_{j=n}^{N}M^{-j}x_j =  \sum_{j=n}^{N}M^{-j}x'_j, \quad n\leq N.
$$
It follows that
$$
  M^Nx =  \sum_{j=0}^{N-n}M^{j}x_{N-j} =  \sum_{j=0}^{N-n}M^{j}x'_{N-j},
$$
and hence
$$
 x_N  - x'_N  = \sum_{j=1}^{N-n}M^{j}(x'_{N-j} - x_{N-j}) = Ms,
$$
where $s\in \mathbb{Z}^d$. This yields $x_N  = x'_N$. Similarly, one can prove   $x_{N-1}  = x'_{N-1}$,  $x_{N-2}  = x'_{N-2}$, etc.\ \ $\Diamond$

\medskip
In what follows we will use the  finite representation for every 
 $x\in X^0$, that is unique by Proposition~\ref{p1.20}. 
 
 It is proved in~\cite{MR09}  that for some matrices $M$, there exist 
 several finite sequences of so-called $M$-nines providing an infinite representation for any  $x\in X^0$. But, as it is said, in this case, we will use only a finite one.
  The situation is more complicated with vectors that have no  finite  expansion~\eqref{1.81}.  In contrast to the positive half-line, for some matrices $M$,  
 there exist  vectors  $x\not\in X^0$ with more than one infinite representations. 
 For example, let  
 $$
  M= \lll \begin{matrix} 1&1\\1&-1  \end{matrix}\rrr,\quad 
 D=\{{\bf 0}, s\}, 
 \quad s=\lll \begin{matrix} 1\\0  \end{matrix}\rrr.
 $$ 
 It is not difficult to see that
 $$
 M^2s=\sum_{k=0}^\infty M^{-2k}s.
 $$
Hence, for every sequence $\{y_k\}_{k=0}^\infty$, where $y_k\in D$, we have
$$
z:=M^2s+\sum_{k=0}^\infty M^{-(2k-1)}y_k=\sum_{k=0}^\infty M^{-2k}s+\sum_{k=0}^\infty M^{-(2k-1)}y_k.
$$ 
Such a vector  $z$ is in $X^0$ only if either the sequence $\{y_k\}_{k=0}^\infty$ is finite, or it contains a  sequence of  $M$-nines (see ~\cite[Theorem~8]{MR09}). It follows that there exists a sequence $\{y_k\}_{k=0}^\infty$ such that the corresponding $z$ is not in 
$X^0$ and it has two different infinite representations~\eqref{1.81}.

\medskip
For every matrix $M\in\mathfrak{M}_d$ and every set of its digits 
 $D$, we set 
$$
U^+ = U^+(M,D) := \{u\in\mathbb{R}^d\,: \ u=\sum_{j=1}^{\infty}M^{-j}u_j, \  u_j\in D \}.
$$
The set $U^+$ is compact, it has a non-empty interior and satisfies the following  self-similarity condition
\begin{equation}
MU^+=\bigcup_{s\in D}(U^+ + s)
\label{1.6}
\end{equation}
(see, e.g., \cite[\S~2.8]{NPS}, \cite[\S~6.1]{KPS}, \cite{B91, LW96}). It is also known that the Lebesgue measure of $U^+$  is an  integer number, but  not necessary equals 1 (methods for computing $\mu(U^+)$ are described in~\cite{GY06}).

\begin{theo} \label{theo1} Representation~\eqref{1.81} is unique for almost all
  $x\in X^+(M,D)$.
\end{theo}

\medskip
{\bf Proof.} According to \eqref{1.6} we have 
$$
   U^+=\bigcup_{s\in D}(M^{-1}U^+ + M^{-1}s);\quad \mu(M^{-1}U^+ + M^{-1}s) = \mu(M^{-1}U^+) = m^{-1} \mu(U^+) 
$$
for all $s\in D$. Hence, the sets $M^{-1}U^+ + M^{-1}s$, $s\in D$, are mutually disjoint up to a set of zero measure. Thus, there exists  a set  $\widetilde{U}\subset U^+$ such that $\mu(\widetilde{U})=0$, and if for some  $u\in U^+\setminus \widetilde{U}$ 
the equality  
$$
   u = \sum_{j=1}^{\infty}M^{-j}u_j= \sum_{j=1}^{\infty}M^{-j}u'_j, \quad  u_j,u'_j\in D,
$$
holds, then $u_1 = u'_1$. Further, let $Mu-u_1\in U^+\setminus \widetilde{U}$. Then  $u\not\in M^{-1}(\widetilde{U}+u_1)$ and 
$$
Mu-u_1 = \sum_{j=1}^{\infty}M^{-j}u_{j+1}= \sum_{j=1}^{\infty}M^{-j}u'_{j+1}, 
$$
which yields $u_2 = u'_2$. Continuing this process and setting 
\begin{equation}
	\widetilde{X}:=\bigcup_{h\in X^0}\bigcup_{n\in\mathbb{Z}}M^n(\widetilde{U} + h)
	\label{1.91},
\end{equation}
 we see that $u_j = u'_j$ for all  $j\ge1$, whenever   $u\not\in  \widetilde{X}$. 
 
Thus we established that any $u\in U\setminus\widetilde{X} $  has a unique representatiom~\eqref{1.81}. 
If now $x\in X^+\setminus \widetilde{X}$, then there exists $k\in\mathbb Z_+$ such that $M^{-k}x\in U$. Since  the relation  $x\not\in  \widetilde{X}$ implies  that  $M^{-k}x\not\in  \widetilde{X}$, and so, as has been proven, the vector $M^{-k}x$  has a unique represantation~\eqref{1.81}, that obviously yields the uniqueness of  $x$. 
It remains to note that $\mu(\widetilde{X})=0$.\quad\quad $\Diamond$   
 
\medskip

\begin{defi}
The space of {\it $M$-positive vectors}, associated with a matrix~$M\in\mathfrak{M}_d$ and a set of its digits $D, $ is defined as follows
$$
   X= X(M,D) := (X^+\setminus \tilde{X})\cup X^0,
$$ 
where  $\widetilde{X}$ is the set given by \eqref{1.91}.
\end{defi}
 Note that every vector $x\in X$ either has a unique representation~\eqref{1.81}, or has a finite representation that will be used whenever we deal with such a vector. 
\medskip
We set
$$
U = U(M,D) := \{x\in X(M,D)\,: \ x_j=\mathbf{0} \quad\mbox{for all}\quad j\le 0 \},
$$
$$
H = H(M,D) := \{x\in  X(M,D)\,: \ x_j=\mathbf{0} \quad\mbox{for all}\quad j>0 \}
$$        
and note that $U$ is an analog of the unit interval $[0,1]$ on the half-line,  and $H$ plays the role  of the set of non-negative integers in  the $M$-adic system. 
The space $X^*$ and the sets \, $U^*$, $H^*$  are defined similarly using the  transposed matrix $M^*$    instead of $M$. For every $\omega\in X^*$ we deal with its  unique 
representation
\begin{equation}
	\omega = \sum_{j=-\infty}^{\infty}(M^*)^{-j}\omega_j, \quad \omega_{j}\in D^*,
	\label{1.92}
\end{equation}
where   only a finite number of nonzero $\omega_{j}$ with negative $j$. Note that  the sets $U$, $U^*$ as well as the whole spaces $X$, $X^*$ have a complex fractal structure for many matrices $M$.

\medskip

	\hspace{-1.5cm}
	\hspace{0.1cm}
	\hspace{.1cm}
	\hspace{0.1cm}
	\hspace{.1cm}

\medskip 
It will be convenient for us to numerate the elements of $H$ and $H^*$ (our "integers"), by the usual non-negative integers. 
For any $k\in\mathbb{Z}_+$ whose  $m$-adic representation is 
$$
k=\sum\limits_{j=0}^{\infty} k_j m^{j}, \quad k_j\in\{0,1,\dots, m-1\},
$$ 
we set
$$
\gamma_{[k]} := \sum\limits_{j=0}^{\infty} M^{j}s_{k_j}, \quad  \gamma^*_{[k]} := \sum\limits_{j=0}^{\infty}  (M^*)^{j}s^*_{k_j}.
$$
Note that \ $H = \{\gamma_{[k]}\,: \ k\in\mathbb{Z}_+\}$ \, and \,
$H^*= \{\gamma^*_{[k]}\,: \ k\in\mathbb{Z}_+\}$, and for every
 $n\in\mathbb{N}$, we set 
$$
H_n := \{\gamma_{[k]}\,: \ 0\leq k\leq m^n-1\}, \quad  H^*_n := \{\gamma^*_{[k]}\,: \ 0\leq k\leq m^n-1\}. 
$$
It is clear that
$$
H_n = \{\gamma\in\mathbb{Z}^d\,: \ \gamma=\sum_{j=1}^{n}M^{j-1}h_j,  \  h_j\in D\} 
$$
and a similar equality holds for $H^*_n$.

Next, we introduce the following sets
\ba
\label{my333}
U_{n,k}:=M^{-n}\gamma_{[k]} + M^{-n}(U),  \quad U_n:= U_{n,0}, \quad  n\in\mathbb{Z},\  k\in\mathbb{Z}_+, 
\ea
and the sets  $U^*_{n,k}$ are defined similarly.

\begin{prop}\label{p1.4}  The sets  $U_{n,k}$,  $U^*_{n,k}$ are almost compact (i.e., are contained in the compact and differ from it by a set of zero measure) and have the following properties:

\medskip 	
(i) For any $n\in\mathbb{N}$ it holds  that
	\begin{equation}
	U = \bigcup_{k=0}^{m^n-1}U_{n,k}, \quad  U^* = \bigcup_{k=0}^{m^n-1}U^*_{n,k}.
	\label{my1.12}
\end{equation}	

\smallskip
(ii) If  $n, n',k, k'\in\mathbb{Z}_+$, $n\ge n'$, then either $U_{n,k}\subset U_{n',k'}$ ($U^*_{n,k}\subset U^*_{n',k'}$), or 
$U_{n,k}$ and $U_{n',k'}$ ($U^*_{n,k}$ and $U^*_{n',k'}$)  are disjoint.
\end{prop}

\noindent
{\bf Proof.} The almost compactness of the set $U$ follows from the compactness
of $U^+$, that obviously implies the almost compactness of any
 $U_{n,k}$. 
 
To prove (i), it suffices to verify that
	\begin{equation}
	M^n(U) = \bigcup_{\gamma\in H_n}(U + \gamma), \quad  (M^*)^n(U^*) = \bigcup_{\gamma^*\in H^*_n}(U^* + \gamma^*).
	\label{1.12}
\end{equation}
For the case  $n=1$, first equality \eqref{1.12}  follows from~\eqref{1.6}, and the inductive step from $n-1$ to $n$ follows from the equality
$$
M^nU = M\bigcup_{\gamma\in H_{n-1}}(U + \gamma) = \bigcup_{\gamma\in H_{n-1}}(MU + M\gamma) = \bigcup_{\gamma\in H_{n-1}}(\bigcup_{s\in D}(U + s) + M\gamma)
$$
$$
= \bigcup_{\gamma\in H_{n-1}}(\bigcup_{s\in D}(U + (s + M\gamma)) = \bigcup_{\gamma\in H_n}(U + \gamma).
$$
To prove (ii) it suffices to consider the case  $n'=0$.
It follows from 1) that there exists  $\gamma_{[l]}\in H$ such that  $U_{n,k}\subset U+ \gamma_{[l]]}$. By the definition of $X$, if   $k'\ne l$, then  the set $U+\gamma_{[k']}$  is disjoint with $U+\gamma_{[l]}$,  hence it is disjoint  with  $U_{n,k}$ as well,  and if, moreover,  $k'= l$, then $U_{n,k}\subset U+k'=U_{n', k'}$.

All statements related to $U^*_{n,k}$ can be proved similarly. 
\ \ $\Diamond$

\medskip
 Obviously,  Proposition~\ref{p1.4} implies the equalities
	\begin{equation}
		X = \bigcup_{\gamma\in H}(U + \gamma), \quad  X^* = \bigcup_{\gamma^*\in H^*}(U^* + \gamma^*).
		\label{1.12'}
	\end{equation}

\medskip
Now, we equip the space $X$ with an  operation of  termwise addition.
The set of digits $D$ of the matrix~$M$ is an abelian group with respect to the addition~$\oplus$ modulo $M$. 
Define the sum of vectors $x,y\in X$ as follows
$$ 
x\oplus y := \sum_{j=-\infty}^{\infty}M^{-j}(x_j\oplus y_j),
$$
where $x_j$,  $y_j$ are digits from  \eqref{1.81} for  $x$ and $y$, respectively, and we choose their finite representation if any. As usual, the equality
 $z=x\ominus y$  means that $z\oplus y = x$. These operations  $\oplus$ и $\ominus$ are similar to the corresponding  operations in the Vilenkin groups and on the half-line
  $\mathbb{R}_+$ (see., e.g., \cite[\S~1.5]{GES}). 
The operation   $\oplus$ on $X^*$  is introduced in a similar way.

\begin{prop}\label{p5.20}
If \, $x,y\in X$, then \, $x\oplus y\in X$.

\end{prop}

\medskip
{\bf Proof.} The statement is trivial if at least one of the vectors $x,y$ or $x\oplus y$ has a finite representation \eqref{1.81}.  Suppose that $x,y$ and $x\oplus y$ do not have finite representations, and let $x'$ be the vector such that $x\oplus x' = \mathbf{0}$. Since $x$ has a unique infinite representation, evidently,  $x'$ also has a unique infinite representation and does not have a finite one. If the vector  $x\oplus y$ has at least two 
infinite representations, then it follows from the equality 
$$
     y = y \oplus (x\oplus x') = (y\oplus x)\oplus x',
$$
that $y$ has at least two infinite representations as well. This contradicts to our assumption.  \ \ $\Diamond$

\bigskip
{\large{\bf \S~3. Walsh functions}
	
	\medskip
	
Given $x\in X$ and  $\omega\in X^*$, we set
\begin{equation}
   \chi(x,\omega) := \exp\left(2\pi i\sum_{j\in {\mathbb Z}}\langle M^{-1} x_j,\,\omega_{1-j}\rangle\right), 
\label{1.02}
\end{equation}
where  $x_j$ and  $\omega_j$  are the digits from decompositions~\eqref{1.81} and~\eqref{1.92}, respectively, and we choose their finite representations if any.
Note that these functions are very similar to the characters for the Vilenkin groups and inherit from them  the following basic  properties: 
$$
  \chi(x\oplus y, \omega) 
= \chi(x,\omega)\cdot\chi(y,\omega), \quad x,y\in X,  \ \omega\in X^*.
$$
$$
\chi(Mx,\omega)=  \chi(x,M^*\omega),  \quad x\in X,  \ \omega\in X^*.
$$

\noindent
These properties follow from \eqref{1.02}  and the fact that
$$
e^{2\pi i \langle M^{-1}(s \oplus t),\, s^*\rangle} =  e^{2\pi i \langle M^{-1}s, \, s^*\rangle} \cdot e^{2\pi i \langle M^{-1}t, \, s^*\rangle}, \quad s,t\in D, \ s^*\in D^*. 
$$

The following statement is related to an analog of the Vilenkin-Chrestenson transform that  is important because of its applications in
function theory and signal processing.

\begin{theo}\label{p6.1}
Let $M\in\mathfrak{M}_d$,  $m=|\det M|$ and $n\in\mathbb{N}$. If  
\begin{equation}
 a_{\gamma^*} = m^{-n} \sum_{\gamma\in H_n} b_{\gamma}\overline{\chi(M^{-n}\gamma,\gamma^*)},  \quad\gamma^*\in H^*_n,
\label{60.1}
\end{equation}
then
\begin{equation}
 b_{\gamma} = m^{-n}  \sum_{\gamma^*\in H^*_n} a_{\gamma^*}\chi(M^{-n}\gamma,\gamma^*),  \quad\gamma\in H_n.
\label{60.2}
\end{equation}
Inversely,  equalities  \eqref{60.1} follow from~\eqref{60.2}.

\end{theo}

\noindent
{\bf Proof.} As above, we use the notation $D=D(M)$ and $D^*=D(M^*)$. Due to Proposition~\ref{p1.2}, the matrix
$$
      \left(\frac{1}{\sqrt{m}}\exp(2\pi i\langle M^{-1}s,s^*\rangle)\right)_{s\in D, s^*\in D^*}
$$ 
is unitary. Taking into account \eqref{1.02}, we see that if  
$$
    a_{s^*} = \sum_{s\in D}b_{s}\overline{\chi(M^{-1}s,s^*)} = \sum_{s\in D}b_{s}\exp(-2\pi i\left\langle M^{-1}s,s^*\rangle)\right),  \quad s^*\in D^*,
$$
then
$$
       b_s = \sum_{s^*\in D^*}a_{s^*}\chi(M^{-1}s,s^*), \quad s\in D.
$$
Thus, the statement is proved for $n=1$. Let us check the inductive step from $n-1$ to $n$. Since every  $x\in H_n$ can be represented in the form $x=\gamma + M^{n-1}s$, where $s\in D$ and  $\gamma\in H_{n-1}$, and a similar representation holds for $\omega\in H^*_n$, for every  $s^*\in D^*$ and $\gamma^*\in H^*_{n-1}$ one can rewrite~\eqref{60.1} as follows
$$
    \alpha_{s^*+M^*\gamma^*}=\frac1{m^{n}}\sum_{\gamma\in H_{n-1}}\sum_{s\in D} \beta_{s+M\gamma}\overline{\chi(M^{-n}(s + M\gamma), \gamma^* + (M^*)^{n-1}s^*)}
$$  
$$
   =\frac1{m^{n}}\sum_{\gamma\in H_{n-1}}\sum_{s\in D} \beta_{s+M\gamma}\overline{\chi(M^{-n}(s \oplus M\gamma), \gamma^* \oplus (M^*)^{n-1}s^*)}.
$$
Hence, using the equalities  $\chi(M^{-n+1}\gamma, (M^*)^{n-1}s^*)=1$, $\chi(M^{-n}s,\gamma^*)=1$ and  ${\chi(M^{-n}s,(M^*)^{n-1}s^*)}=\chi(M^{-1}s, s^*),$   we get
$$
    \alpha_{s^*+M^*\gamma^*} = \frac{1}{m^{n-1}}\sum_{\gamma\in H_{n-1}}\left(\frac{1}{m}\sum_{s\in D}\beta_{s+M\gamma}\overline{\chi(M^{-n}s,(M^*)^{n-1}s^*)}\right)\overline{\chi(M^{-n+1}\gamma,\gamma^*)},
$$
that, due to the induction hypothesis, yields
$$у
   \frac{1}{m}\sum_{s\in D}\beta_{s+M\gamma}\overline{\chi(M^{-1}s, s^*)}=\frac{1}{m^{n-1}}\sum_{\gamma^*\in H^*_{n-1}}\alpha_{s^*+M^*\gamma^*}\chi(M^{-n+1}\gamma,\gamma^*).
  $$
  Next, using the induction base and the same equalities, 
  we obtain
$$
  \beta_{s+M\gamma} = \frac1{m^{n}}\sum_{s^*\in D^*}\sum_{\gamma^*\in H^*_{n-1}} \alpha_{s^*+M^*\gamma^*}\chi(M^{-n+1}\gamma,\gamma^*)\chi(M^{-1}s,s^*) 
$$
$$
   = \frac1{m^{n}}\sum_{s^*\in D^*}\sum_{\gamma^*\in H^*_{n-1}} \alpha_{s^*+M^*\gamma^*}\chi(M^{-n}(s + M\gamma),\gamma^* + (M^*)^{n-1}s^*).
$$
Similarly, using~\eqref{60.2}, one can  prove~\eqref{60.1}.
\ \ \ $\Diamond$

\medskip
\begin{defi}
The {\it Walsh system} is the set $\{W_{\alpha}\}_{\alpha=0}^{\infty}$, where the functions $W_{\alpha}$ are defined on $X$ by the formula
$$
      W_{\alpha}(x) :=  \chi(x, \gamma^*_{[\alpha]}),  \quad \alpha\in\mathbb{Z}_+,  \ x\in X
$$
\end{defi}

\begin{theo} \label{theo3.1} The Walsh function $W_{\alpha}$ \textcolor{blue} is $H$-periodic; moreover,  if $0\leq \alpha, k \le m^{n}-1$, then  $W_{\alpha}$ is constant on $U_{n,k}$.

\end{theo}


\medskip
\noindent
{\bf Proof.} 
Let $ x = u + h, \  u\in U, \  h\in H$. Then, since  $u+h=u\oplus h$ and
$\chi(h, \gamma^*_{[\alpha]})=1$, we have
$$
W_{\alpha}(x)=\chi(u+h, \gamma^*_{[\alpha]})=\chi(u, \gamma^*_{[\alpha]})\chi(h, \gamma^*_{[\alpha]})=W_\alpha(u),
$$
that just means the $H$-periodicity of the function $W_{\alpha}$.

 Now, let $0\le \alpha, k\le m^{n}-1$ and $x\in U_{n,k}$, i.e. 
$x=M^{-n}\gamma_{[k]}+M^{-n}u$, $u\in U$. 
 Then, using the equalities  $M^{-n}\gamma_{[k]}+M^{-n}u=M^{-n}\gamma_{[k]}\oplus M^{-n}u$
and $\chi(M^{-n}u, \gamma^*_{[\alpha]})=1$, we obtain
 $$
   W_{\alpha}(x)=\chi(x, \gamma^*_{[\alpha]})=\chi(M^{-n}u+M^{-n}\gamma_{[k]}, \gamma^*_{[\alpha]}) = \chi(M^{-n}\gamma_{[k]},\gamma^*_{[\alpha]}),
$$ 
that means that the function $W_{\alpha}$ is constant on the set $U_{n,k}$.
    \ \ \ \
$\Diamond$

\begin{prop}\label{p6.2}
	The Walsh system $\{W_{\alpha}\}_{\alpha=0}^{\infty}$ is orthogonal  in $L^2(U)$, i.e.,
	$$
	\frac{1}{\mu(U)} \int\limits_U W_{\alpha}(x)\overline{W_{\beta}(x)}\,d\mu(x) = \delta_{\alpha,\beta},  \quad \alpha,\beta\in\mathbb{Z}_+,
	$$
	where  $\delta_{\alpha,\beta}\, $ is the Kronecker delta. 
\end{prop}

This statement is the direct corollary of the following lemma.

\begin{lem}\label{lm1.1}
{If $\omega\in H^*$, $\omega\neq\mathbf{0}$, then} 
\begin{equation}
 \int\limits_U\chi(x,\omega)\,d\mu(x) = 0.
\label{1.14}
\end{equation}
\end{lem}

\noindent
{\bf Proof.} Let $\omega\in H^*$, $\omega\neq\mathbf{0}$. Suppose first that 
 $\omega$ is such that its digit $\omega_0$ in expansion~\eqref{1.91} is nonzero. Then, using~\eqref{1.12} and taking into account that $M^{-1}x+M^{-1} s=M^{-1}x\oplus M^{-1} s$ whenever  $x\in  U$ $s\in D$, we obtain 
$$
  \int\limits_U\chi(x,\omega)\,d\mu(x) = \frac1m\int\limits_{MU}\chi(M^{-1}x,\omega)\,d\mu(x) = \frac1m\sum_{s\in D}\int\limits_{U+s}\chi(M^{-1}x,\omega)\,d\mu(x) 
$$
$$
   = \frac1m\sum_{s\in D}\int\limits_{U}\chi(M^{-1}x + M^{-1}s),\omega)\,d\mu(x) = \frac1m\sum_{s\in D}\chi(M^{-1}s,\omega)\int\limits_{U}\chi(M^{-1}x,\omega)\,d\mu(x). 
$$
By Proposition~1, it follows that
$$
\sum_{s\in D}\chi(M^{-1}s,\omega) = \sum_{s\in D}\exp(2\pi i \left\langle M^{-1}s,\omega_{0} \right\rangle) = 0.
$$
Thus,  equality~\eqref{1.14} is proved for the case $\omega_0\neq \mathbf{0}$.

Now, we suppose that $\omega_0 = \omega_{-1}=\dots =\omega_{-n}=\mathbf{0}$ and   $\omega_{-n}\neq\mathbf{0}$. 
Since $\chi(\gamma,\omega)=1$  for all $\gamma\in H$, it follows from Proposition~\ref{p1.4} that
$$
  \int\limits_{M^nU}\chi(x,\omega)\,d\mu(x) = \sum_{\gamma\in H_n}\int\limits_{U+\gamma}\chi(x,\omega)\,d\mu(x) = \sum_{\gamma\in H_n}\int\limits_{U}\chi(x + \gamma,\omega)\,d\mu(x) 
$$
$$
 = \sum_{\gamma\in H_n}\int\limits_{U}\chi(x \oplus \gamma,\omega)\,d\mu(x) 
= \sum_{\gamma\in H_n}\int\limits_{U}\chi(x,\omega)\chi(\gamma,\omega)\,d\mu(x) = m^n  \int\limits_{U}\chi(x,\omega)\,d\mu(x).
$$
So, to prove~\eqref{1.14} it remains to verify that $$\int\limits_{M^nU}\chi(x,\omega)\,d\mu(x) = 0.$$
Since
${(M^*)}^{-n}\omega\in H^*$  and  ${({(M^*)}^{-n}\omega)}_0\neq\mathbf{0}$, using the above conclusion, we have  $\int\limits_{U}\chi(x,{(M^*)}^{-n}\omega)\,d\mu(x) = 0$. It remains to note that
$$
  \int\limits_{U}\chi(x,{(M^*)}^{-n}\omega)\,d\mu(x) = \int\limits_{U}\chi(M^{-n}x,\omega)\,d\mu(x) = m^n \int\limits_{M^nU}\chi(x,\omega)\,d\mu(x).\ \ \ \Diamond
$$

\begin{lem}\label{p1.7}
For every $n\in\mathbb{N}$, the following equality holds  
\begin{equation}
m^{-n}\sum_{\gamma^*\in H^*_n}\chi(x,\gamma^*)=\left\{
\begin{array}{l}
\ 1,\quad \ x\in U_n,\\
\ 0, \quad \ x\in U\setminus U_n.
\end{array}
\right.
\label{1.16}
\end{equation}

\end{lem}

\noindent
{\bf Proof.} If $x\in U_n$ and $\gamma^*\in  H^*_n$, then  $\chi(x,\gamma^*) = 1$,   
 which yields~\eqref{1.16} whenever  $x\in U_n$. 
 
 Suppose that $\gamma^*\in H^*_n$, $x\in U\setminus U_n$, 
and $x_{j_0}\neq\mathbf{0}$ for some  $j_0\le n$. Using the equalities $x=M^{-1}x_1\oplus M^{-2}x_2\oplus\dots \oplus M^{-n}x_n\oplus x'$ and
$\chi(x',\gamma^*) = 1$,  we get
$$
\sum_{\gamma^*\in H^*_n}\chi(x,\gamma^*)= \sum_{\gamma^*\in H^*_n}\chi(M^{-1}x_1,\gamma^*)\dots \chi(M^{-n}x_n,\gamma^*)
$$
$$
= \sum_{\gamma^*\in H^*_n}\chi(M^{-1}x_1,\gamma^*_0)\dots \chi(M^{-n}x_n, (M^*)^{n-1}\gamma^*_{1-n})
$$
$$
= \sum_{\gamma_1^*\in D^*}\chi(M^{-1}x_1,\gamma^*_0)\dots \sum_{\gamma^*_{1-n}\in D^*}\chi(M^{-n}x_n, (M^*)^{n-1}\gamma^*_{1-n}). 
$$
To prove~\eqref{1.16}, it remains to note that  Proposition~1 implies
$$
    \sum_{\gamma_{1-j_0}^*\in D^*}\chi(M^{-j_0}x_{j_0},(M^*)^{j_0-1}\gamma^*_{1-j_0})= \sum_{\gamma_{1-j_0}^*\in D^*}\exp\left(2\pi i\langle M^{-1}x_{j_0},\gamma^*_{1-j_0}\rangle\right) = 0.\ \ \ \Diamond
$$

\medskip 
{\bf Remark.}
Suppose that $x\in U_{n,k}$ for some $k\le m^n-1$, i.e.,
$$
x = M^{-n}\gamma_{[k]} + x', \quad x'=M^{-n}y,
$$
where $y\in U$, and let $\gamma^*\in H^*_n$.  Then
$$
\chi(x,\gamma^*)= \chi(M^{-n}(\gamma_{[k]} \oplus y),\gamma^*)= \chi(M^{-n}\gamma_{[k]},\gamma^*)\chi(x',\gamma^*).
$$
Since $x'\in U_n$, as a consequence of Lemma \ref{p1.7} we get
\begin{equation}
	\sum_{\gamma^*\in H^*_n}c_k(\gamma^*)\chi(x,\gamma^*)=\left\{
	\begin{array}{l}
		\ 1,\quad \ x\in U_{n,k},\\
		\ 0, \quad \ x\in U\setminus U_{n,k},
	\end{array}
	\right.
	\label{11.1}
\end{equation}
where $c_k(\gamma^*)=\chi(M^{-n}\gamma_{[k]},\gamma^*)$, $0\leq k\leq m^n-1$.

\begin{theo} \label{theo2} 
	The Walsh system  $\{W_{\alpha}\}_{\alpha=0}^{\infty}$ 
is an orthogonal basis of the space	 $L^2(U)$.
\end{theo}

\medskip
\noindent
{\bf Proof.} 
Due to Proposition~\ref{p6.2}, it suffices to verify that the system   $\{W_{\alpha}\}_{\alpha=0}^{\infty}$ is complete.
Suppose that a function $f\in L^2(U)$ satisfies the condition
\begin{equation}
	\int\limits_Uf(x)\overline{\chi(x,\gamma^*)}\,d\mu(x) = 0, \quad \gamma^*\in H^*
	\label{11.2}
\end{equation}
and prove, that $f$ equals zero almost everywhere  (we set $f(x)=0$ for $x\in X\setminus U$).
By \eqref{11.1} and \eqref{11.2}, for $n\in \mathbb{N}$ and $k\in\{0,1,\dots,m^n-1\}$ we have 
\begin{equation}
	\int\limits_{U_{n,k}}f(x)\,d\mu(x) = \sum_{\gamma^*\in H^*_n}c_k(\gamma^*)\int\limits_Uf(x)\overline{\chi(x,\gamma^*)}\,d\mu(x) = 0.
	\label{11.3}
\end{equation}

Now, let $G$ be an open ball in $\mathbb{R}^d$, and let  
 $n_1$ denote the minimal number~$n$ for which there exists $k$ such that $U_{{n},k}\subset G\cap U$. Denote by  $\Delta_1,\dots, \Delta_{s_1}$ all such sets $U_{{n_1},k}$.  Let us denote by
 $n_2$  the minimal number $n > n_1$ for which there exists a set $U_{n,k}$ contained in  $G\cap U$, but not contained in any set  $\Delta_l$, $l\le s_1$. Let  $\Delta_{s_1+1},\dots, \Delta_{s_2}$ 
 be all such sets. Continuing this process, we obtain a sequence of "cells" \ $\Delta_j$ that, by  Proposition~\ref{p1.4}, are pairwise disjoint. It follows from the construction that the union of all such  $\Delta_j$ is contained in $G\cap U$. Let us check the inverse.
Let $x_0\in G\cap U$. It follows from Proposition~\ref{p1.4} that  for every $n$ there exists $k$ such that $x_0\in U_{{n},k}$. Due to the boundedness of $U$ and \eqref{my000},  the set   $M^{-n}(U)$ is contained in
an arbitrarily small neighborhood of zero whenever~$n$ is large enough. Therefore, a set 
$U_{n,k}$ containing a point $x_0$  is contained in an arbitrarily small neighborhood 
of $x_0$ whenever~$n$ is large enough.
 Since $G$ is open, it follows that some sets $U_{n,k}$ are contained in $ G\cap U$ 
 together with a point $x_0\in  G\cap U$. Choosing from all such sets
the one that corresponds to the smallest $n$, we see that it is just one of the "cells" of $\Delta_j$.  Thus, we have
$$
\bigcup_{j=1}^{\infty}\Delta_j= G\cap U. 
$$
 Since  $U_{n,k}$ is contained in  $U$ if and only if  $0\leq k\le m^n-1$,  it follows from~\eqref{11.3} that 
$$
\int\limits_{G}f(x)\,d\mu(x) = 0.
$$ 
Hence, if now $x^*$ is a Lebesgue point of  $f$ and   $G_r(x^*)$ is an open ball of radius~$r$ with a center  $x^*$, then
$$
\lim\limits_{r\to0}\frac{1}{\mu G_r(x^*)}\int\limits_{G_r(x^*)}|f(x)-f(x^*))|\,d\mu(x)=0,\quad \int\limits_{G_r(x^*)}f(x)\,d\mu(x) = 0,
$$
which yields $f(x^*)=0$.  It remains to note that almost all points are the 
 Lebesgue points of a localy summable function $f$. 
\ \ $\Diamond$

\bigskip

\medskip
A finite linear combination of the Walsh functions is called a 
{\it Walsh polynomial} (on $X$). 
We will say that 
$$
                w(x) = \sum\limits_{k=0}^{m^{n} - 1}a_{k}W_k(x), \quad  x\in X,
$$
is a Walsh polynomial  of order $n$ if
 $a_{j}\neq 0$ for some $j\in\{m^{n-1},\dots,m^n-1\}$. The Walsh polynomials on  $X^*$ are defined in a similar way.

\begin{prop}\label{p11.2}
A function $w :\, X\to\mathbb{C} $ is a Walsh polynomial of order $n$ if and only if 
it is an $H$-periodic function that is constant on each set  $U_{n,k}$,  $0\leq k \leq m^n-1$. 
\end{prop}

This statement follows from Theorem~\ref{theo3.1} and Proposition~\ref{p6.1}.

\bigskip

{\large{\bf \S~4. Fourier transform}
	
\begin{defi}
 For a function 	  $f\in L^1(X)$, its {\it Fourier transform} is defined by 
\begin{equation}
   \widehat f(\omega)=\frac{1}{\mu(U)}\int\limits_Xf(x)\overline{\chi (x,\omega)}d\mu(x), \quad \omega\in X^*.
\label{1.03}
\end{equation}
\end{defi}
 Note that a similar concept on  $\mathbb{R}_+$  is known \cite{SWS} as  the 
Walsh-Fourier transform.

\medskip
Let us introduce the following classes of functions  $f$ given on $X$. 
$$
 \mathcal{S}_n(X) := \{ f\,: \ f(x)=f(x')\quad\mbox{for all}\quad x, x'\in U_{n,k},
  \  k\in {\mathbb Z}_+\},
$$
$$
   \mathcal{S}_n^{(k)}(X) := \{ f\in \mathcal{S}_n(X)\,: \ \supp f \subset M^k(U) \}, \quad
       \mathcal{S}(X) := \bigcup_{n,k} \mathcal{S}_n^{(k)}(X).
$$
The classes $\mathcal{S}_n(X^*)$,  $\mathcal{S}_n^{(k)}(X^*)$ and $\mathcal{S}(X^*)$  are defined similarly.

Note that the elements of $\mathcal{S}_n^{(k)}(X)$ are step functions, and $k$ defines the size of the step.
Also, any function $f\in \mathcal{S}_n(X)$ can be represented in the form
\begin{equation}
	f =\sum\limits_{k=0}^{\infty}f_{n,\, k}{\bf 1}_{U_{n,\,k}}, 
	\label{2.18}
\end{equation}
where $f_{n,\,k}:= f(M^{-n}\gamma_{[k]})$  is the value of  $f$ on the set  $U_{n,\,k}$.

\medskip
\begin{prop}\label{p7.2} 
	The following properties hold:
	
\smallskip
\textup{(a)}  if \,$f\in L^1(X)\cap {\mathcal S}_n(X),$ then $\supp\widehat f \subset (M^*)^{n}U^*$\textup{;}  

\smallskip
\textup{(b)}  if \, $f\in L^1(X)$\, and\, $\supp f \subset M^nU$,  then $\widehat f \in {\mathcal S}_n(X^*)$.
\newline

\end{prop}

\smallskip
{\bf Proof.} Let $f\in L^1(X)\cap {\mathcal S}_n(X)$. Then, due to~\eqref{2.18},
for every  $\omega\in X^*$, there holds 
$$
    \widehat f(\omega) =\frac{1}{\mu(U)}  \int\limits_X f(t)\overline{\chi(t,\omega)}\,d\mu(t) = \frac{1}{\mu(U)} \sum_{k=0}^{\infty}f_{n,\,k}\int\limits_{U_{n,\,k}}\overline{\chi(t,\omega)}\,d\mu(t).
$$
If $t\in U_{n,\,k}$, then $t=M^{-n}\gamma_{[k]}+M^{-n}x$,  where $x\in U$,   and taking into account that $M^{-n}\gamma_{[k]}+M^{-n}x=M^{-n}\gamma_{[k]}\oplus M^{-n}x$, we obtain 
$$
    \int\limits_{U_{n,\,k}}\chi(t,\omega)\,d\mu(t) =  m^{-n}\int\limits_{U}\chi(M^{-n}\gamma_{[k]}+M^{-n}x,\omega)\,d\mu(x)
    $$
    $$
   =\frac{\chi(M^{-n}\gamma_{[k]},\omega)}{m^n}\int\limits_{U}\chi(M^{-n}x,\omega)\,d\mu(x) =\frac{\chi(M^{-n}\gamma_{[k]},\omega)}{m^n}\int\limits_{U}\chi(x,\lambda)\,d\mu(x),
$$
where $\lambda=(M^*)^{-n}\omega$. If now 
$\omega\not\in(M^*)^{n}U^*$, then $\lambda=h+y$, $ y\in U^*$, $h\in H^*$, where  $h\ne\mathbf{0}$, and taking into account that  $h+y=x\oplus y$ and $\chi(x, y)=1$, 
it follows from Lemma~\ref{lm1.1} that
$$
\int\limits_{U}\chi(x,\lambda)\,d\mu(x)=\int\limits_{U}\chi(x,h\oplus y)\,d\mu(x)=\int\limits_{U}\chi(x,h)\,d\mu(x)=0.
$$
Comparing this equality with the two previous ones, we get $\widehat f(\omega) = 0$, and thus assertion (a) is proved.

Suppose now that  $f\in L^1(X)$, $\supp f \subset M^nU$, and let $\omega\in U^*_{n,\,k}$, i.e.,
 $\omega=(M^*)^{-n}(\gamma^*_{[k]}+y)$,  where $y\in U^*$. 
 Then, taking into account that  $\gamma^*_{[k]}+y=\gamma^*_{[k]}\oplus y$, 
 $\chi(M^nx, \omega)= \chi(x, (M^*)^n\omega)$ and
 $\chi(x,y)=1$,  for every $x\in U$, we obtain 
$$
    \widehat f(\omega) = \frac{1}{\mu(U)} \int\limits_{M^nU} f(t)\overline{\chi(t,\omega)}\,d\mu(t) = 
     \frac{m^n}{\mu(U)} \int\limits_{U} f(M^nx)\overline{\chi(M^nx,\omega)}\,d\mu(x)
$$
$$
 = \frac{m^n}{\mu(U)} \int\limits_{U} f(M^nx)\overline{\chi(x,\gamma^*_{[k]}+y)}\,d\mu(x)=
\frac{m^n}{\mu(U)} \int\limits_{U} f(M^nx)\overline{\chi(x,\gamma^*_{[k]})}\,d\mu(x)
$$
$$
=\frac{m^n}{\mu(U)} \int\limits_{U} f(M^nx)\overline{\chi(M^n x,(M^*)^{-n}\gamma^*_{[k]})}\,d\mu(x)
= \widehat f((M^*)^{-n}\gamma^*_{[k]}),
$$
which means that $ \widehat f$ is constant on the set $ U^*_{n,\,k}$, i.e.,  $\widehat f \in {\mathcal S}_n(X^*)$.
\ \ $\Diamond$

\begin{theo}\label{p8.2}
For every integers  $n$ and $k$, there holds 
$$
    f \in \mathcal{S}_n^{(k)}(X) \ \Longleftrightarrow \ \widehat{f} \in \mathcal{S}_k^{(n)}(X^*). 
$$
\end{theo}
This statement is a direct consequence of Proposition~\ref{p7.2}.

\begin{prop}\label{p99.2}
The set  $\mathcal{S}(X)$  is dense in $L^2(X)$.
\end{prop}

\smallskip
{\bf Proof.}
Let  $f\in L^2(X)$ и $\varepsilon>0$. Because  of the absolute continuity of the Lebesgue integral, there exists an integer $n$ such that
\begin{equation}
    \int\limits_{X\setminus M^n(U)}|f(x)|^2\,d\mu(x) < \varepsilon. 
\label{17}
\end{equation}
Setting $\widetilde f := f|_{M^n(U)}$, we have $\supp \widetilde f(M^n\cdot)\subset U$ and  $\widetilde f(M^n\cdot)\in L^2(U)$.
Due to Theorem~\ref{theo2}, there exist coefficients  $c_0, c_1,\dots, c_N$ such that
$$
   \int\limits_U| \widetilde f(M^nx) - \sum_{\alpha=0}^Nc_{\alpha}W_{\alpha}(x)|^2\,d\mu(x) < \frac{\varepsilon}{m^n}.
$$
Hence, if now
$$
                 g(x):= \sum_{\alpha=0}^Nc_{\alpha}W_{\alpha}(M^{-n}x)\Big|_{M^n(U)},
$$
then
\begin{equation}
    \int\limits_{M^n(U)}|f(x) - g(x)|^2\,d\mu(x) < \varepsilon. 
\label{18}
\end{equation}
Since   $g\in\mathcal{S}(X)$, the density of $\mathcal{S}(X)$ in $L^2(X)$ 
follows from~ \eqref{17},  \eqref{18} and Proposition~\ref{p1.4}. \ \  $\Diamond$

\smallskip
{\bf Remark.} It is clear from the proof of Proposition~\ref{p99.2} 
that any compactly supported function  $f\in L^2(X)$ can be approximated by functions from the space   $\mathcal{S}(X)$ with the same support.
 
\medskip
The next statement is the Poisson summation formula for the functions from
 $\mathcal{S}(X)$; for some other generalizations of the  classical Poisson summation formula, see, e.g., \cite{Pa11}.

\begin{prop}\label{p95}
For any  $f \in \mathcal{S}(X)$, the following equality   
\begin{equation}
    \sum_{\gamma\in H}f(\gamma) =  \sum_{\gamma^*\in H^*}\widehat{f}(\gamma^*) 
\label{95.1}
\end{equation}
holds true.

\end{prop}

\smallskip
{\bf Proof.} The function
\begin{equation}
    g(x) := \sum_{\gamma\in H}f(x + \gamma), \quad x\in X, 
\label{95.2}
\end{equation}
is an  $H$-periodic step-function. Hence, by Theorem~\ref{theo2}, it can be  decomposed as follows
\begin{equation}
    g(x) = \sum_{\gamma^*\in H^*}a_{\gamma^*}(g)\chi(x,\gamma^*), \quad x\in U, 
\label{95.3}
\end{equation}
where 
$$
  a_{\gamma^*}(g) = \frac{1}{\mu(U)}\int\limits_U g(x)\overline{\chi(x,\gamma^*)}\, d\mu(x).
$$ 
Applying  \eqref{1.12'}  and \eqref{95.2}, we obtain 
$$
   a_{\gamma^*}(g) = \frac{1}{\mu(U)}\sum_{\gamma\in H}\int\limits_U f(x + \gamma)\overline{\chi(x,\gamma^*)}\, d\mu(x)
$$
$$
    = \frac{1}{\mu(U)}\sum_{\gamma\in H}\int\limits_{U + \gamma} f(x)\overline{\chi(x,\gamma^*)}\, d\mu(x)
$$
$$
   = \frac{1}{\mu(U)}\int\limits_Xf(x)\overline{\chi (x,\gamma^*)}d\mu(x) = \widehat{f}(\gamma^*).
$$
This together with~\eqref{95.3}, where $x=\mathbf{0}$, yields \eqref{95.1}. \  $\Diamond$

\begin{prop}\label{p9.2}
Let  $f \in \mathcal{S}(X)$. Then for every $x\in X$, there holds
$$
    f(x) = \frac{1}{\mu(U^*)} \int\limits_{X^*} \widehat{f}(\omega)\chi(x,\omega)\,d\mu(\omega).
$$

\end{prop}

\smallskip
{\bf Proof.} Due to \eqref{2.18}, it suffices to consider the case where  $f$ is the characteristic function of a set $U_{n,k}$,  $n\in {\mathbb Z}$,  
$k\in {\mathbb Z}_+$.
In this case, for every  $x\in X$, using Theorem~\ref{p8.2}, we have
$$
   I(x):= \frac{1}{\mu(U^*)} \int\limits_{X^*} \widehat{f}(\omega)\chi(x,\omega)\,d\mu(\omega)
$$   
$$   
   = \frac{1}{\mu(U^*)}\int\limits_{(M^*)^nU^*}\,d\mu(\omega)\left(\frac{1}{\mu(U)}\int\limits_{U_{n,\,k}} \overline{\chi(t,\omega)}\,d\mu(t)\,\chi(x,\omega)\right) 
$$
$$
    = \frac{1}{\mu(U^*)}\int\limits_{U^*}\,d\mu(\omega)\left(\frac{1}{\mu(U)}\int\limits_{U} \overline{\chi(M^{-n}(t\oplus \gamma_{[k]}),(M^*)^n\omega)}\,d\mu(t)\,\chi(x,(M^*)^n\omega)\right) 
$$
$$
    = \frac{1}{\mu(U^*)}\int\limits_{U^*}\,d\mu(\omega)\left(\frac{1}{\mu(U)}\int\limits_{U} \overline{\chi(t,\omega)}\,d\mu(t)\,\chi(M^nx \ominus\gamma_{[k]},\omega)\right). 
$$
Since  $\chi(t,\omega)=1$  whenever $t\in U$ $\omega\in U^*$, it follows that
$$
  I(x) =  \frac{1}{\mu(U^*)}\int\limits_{U^*}\chi(M^nx \ominus \gamma_{[k]},\omega)\,d\mu(\omega).
$$
Thus, if $x\in U_{n,\,k}$, then $M^nx \ominus \gamma_{[k]}\in U$ and $I(x)=1$. In the case of  $x\notin U_{n,\,k}$, we have 
$M^nx \ominus \gamma_{[k]} = x' + x''$, where $x'\in U$, $x''\in H$, $x''\neq \mathbf{0}$, and since $\chi(x',\omega)=1$, $x'+x''=x'\oplus x''$, by Lemma~\ref{lm1.1} it follows that
$$
 I(x) =  \frac{1}{\mu(U^*)}\int\limits_{U^*}\chi(x',\omega)\chi(x'',\omega)\,d\mu(\omega)=0.
$$
Thus, we established the equality  $I=f$. \quad  $\Diamond$

\medskip
The transform $f\to\widehat{f}$, given on $L^1(X)\cap L^2(X)$ by formula~\eqref{1.03}, can be extended to the space 
 $L^2(X)$ in a standard way. 
Namely, if  $f\in L^2(X)$ and 
\begin{equation}
    J_lf(\omega):= \frac{1}{\mu(U)}\int\limits_{M^l(U)}f(x)\overline{\chi(x,\omega)}\,d\mu(x),  \quad l\in\mathbb{Z}_+, \ \omega\in X^*,
\label{2.17}
\end{equation}
then $\widehat{f}$ coincides with  the limit of $J_lf$  in $L^2(X^*)$ as $l\to +\infty$.

The following analog of the Plancherel theorem holds.

\begin{theo} \label{theo2.2}
If $f\in L^2(X)$, then  $\widehat{f}\in L^2(X^*)$ and  
\begin{equation}
   \frac{1}{\mu(U)}\int\limits_{X}|f(x)|^2\,d\mu(x) = \frac{1}{\mu(U^*)}\int\limits_{X^*}|\widehat{f}(\omega)|^2\,d\mu(\omega).
\label{12.1}
\end{equation}
If, moreover, $f,g \in L^2(X)$, then  we have
\begin{equation}
\frac{1}{\mu(U)}\int\limits_{X}f(x)\overline{g(x)}\,d\mu(x) = \frac{1}{\mu(U^*)}\int\limits_{X^*}\widehat{f}(\omega)\overline{\widehat{g}(\omega)}\,d\mu(\omega).
\label{12.2}
\end{equation}

\end{theo}

\smallskip
{\bf Proof.} Let first $f\in\mathcal{S}(X)$. Then, by Definition~\ref{p9.2}, we obtain 
$$
\frac{1}{\mu(U)}\int\limits_{X}f(x)\overline{f(x)}\,d\mu(x) = \frac{1}{\mu(U)}\int\limits_{X}f(x)\,d\mu(x) \frac{1}{\mu(U^*)} \int\limits_{X^*} \overline{\widehat{f}(\omega)}\overline{\chi(x,\omega)}\,d\mu(\omega)
$$
$$
= \frac{1}{\mu(U^*)} \int\limits_{X^*} \overline{\widehat{f}(\omega)}\,d\mu(\omega) \frac{1}{\mu(U)}\int\limits\limits_{X}f(x)\overline{\chi(x,\omega)}\,d\mu(x)= \frac{1}{\mu(U^*)}\int\limits_{X^*}|\widehat{f}(\omega)|^2\,d\mu(\omega).
$$
Thus,  equality~\eqref{12.1} is proved for the case of $f\in\mathcal{S}(X)$. 

\smallskip
Let now $f\in L^2(X)$ and $\supp f\subset M^l(U)$ for some $l\in\mathbb{Z}$.  Because of  Proposition~\ref{p99.2} and the remark to it, there exists a sequence  $\{f_k\}_{k=0}^{\infty}$ converging to $f$ in the norm of $L^2(X)$ and such that  $f_k\in\mathcal{S}(X)$  and  $\supp f_k\subset M^l(U)$. For any  $k,n\in\mathbb{Z}_+$, we have 
$$
  \frac{1}{\mu(U)}\int\limits_{X}|f_k(x)-f_n(x)|^2\,d\mu(x) = \frac{1}{\mu(U^*)}\int\limits_{X^*}|\widehat{f_k}(\omega) - \widehat{f_n}(\omega)|^2\,d\mu(\omega)
$$
(here relation~\eqref{12.1} was applied to $f_k-f_n\in\mathcal{S}(X)$). 
Hence the sequence $\{\widehat{f}_k\}_{k=0}^{\infty}$ is fundamental and converges in $L^2(X^*)$ 
to some function $F$. In addition, since
 $\supp f\subset M^l(U)$,  $\supp f_k\subset M^l(U)$ and 
 $f_k\overset {L^2}\to f$, the sequence $\{\widehat{f}_k\}$ converges to $\widehat{f}$ uniformly on any compact set $K\subset X^*$. This yields that the functions $F$ and $\widehat f$ are equivalent. Therefore, one can obtain~\eqref{12.1}  by passing to the limit as $k\to\infty$ from the equalities
 $$
     \frac{1}{\mu(U)}\int\limits_{X}|f_k(x)|^2\,d\mu(x) = \frac{1}{\mu(U^*)}\int\limits_{X^*}|\widehat{f_k}(\omega)|^2\,d\mu(\omega), \quad k\in\mathbb{Z}_+.
$$
Now, we consider an arbitrary $f\in L^2(X)$, and let $f_l$ be the restriction of $f$ onto $M^l(U)$. Since the theorem is proved already for
 $f_l$,  passing to the limit as $l\to\infty$ in the equality 
$$
     \frac{1}{\mu(U)}\int\limits_{X}|f_l(x)|^2\,d\mu(x) = \frac{1}{\mu(U^*)}\int\limits_{X^*}|\widehat{f_l}(\omega)|^2\,d\mu(\omega),
$$
we get~\eqref{12.1}.
 It remains to note that equality~\eqref{12.2} can be deduced from~\eqref{12.1} exactly as for the Fourier transform in   $L^2(\mathbb{R}^d)$ (see, e.g., \cite[Theorem 8.15]{Igari}).~$\Diamond$
 
 \medskip

 \begin{theo}\label{p71.200}
 	For every  $g\in L^2(X^*)$ there exists $f\in L^2(X)$ such that $ \widehat f=g$.  
 \end{theo}
 
 {\bf Proof.} First we consider the case where $g\in {\mathcal S}(X^*)$, and let 
 $$
 \check g(x): = \frac{1}{\mu(U^*)} \int\limits_{X^*} g(\omega)\chi(x,\omega)\,d\mu(\omega).
 $$
 Similarly  to Proposition~\ref{p9.2} (dealing with $\check g$ instead of $\hat f$), we obtain  
 $$
 g(\omega)=\frac{1}{\mu(U)} \int\limits_{X}\check g(x)\overline{\chi(x,\omega)}\,d\mu(x),
 $$
 which means that  $f=\check g$ is a required function $f$, and the mapping $g\to \check g$ can be naturally called by the inverse Fourier transform. 
 
 Let now $g$ be an arbitrary function from $L^2(X^*)$. 
 Because of Proposition~\ref{p99.2} for $X^*$ instead of $X$, 
 given $\varepsilon>0$ there exists $\tilde g\in {\mathcal S}(X^*)$ such that $\|g-\tilde g\|<\varepsilon$. The inverse Fourier transform can be extended onto $L^2(X^*)$ in the same way as it was done for the Fourier transform. 
 Note that changing roles of $X$ and $X^*$, an analog of Theorem~\ref{p9.2} for  $g$ and $\check g$ instead of $f$ and  $\hat f$ can be easily  checked. Because of this, we have $\|\check g-\check {\tilde g}\|<\frac{\mu U}{\mu U^*}\varepsilon$. Next, it follows from 
 Theorem~\ref{p9.2} that $\|\widehat{\check g}-\widehat{\check {\tilde g}}\|<\varepsilon$.  Setting $f:= \check g$ and using that the equality  $\widehat{\check{\tilde g}}=\tilde g$ is already proved, we obtain $\|\widehat f-g\|<2\varepsilon$, which yields $ \widehat f=g$. \ \ $\Diamond$
 
 \

\medskip

\begin{theo}\label{p71.2}
Let  $f\in L^2(X)$. For $ {f\in\mathcal S}_n(X)$, it is necessary and sufficient that $\supp\widehat f \subset (M^*)^n(U^*)$.  
\end{theo}

\smallskip
{\bf Proof.} {\it Necessity.} Let $f\in L^2(X)\cap {\mathcal S}_n(X)$, and let   $f_k$ be a restriction of $f$ onto the set $M^k(U)$.   Since $f_k\in L^1(X)\cap {\mathcal S}_n(X)$, it follows from Proposition~\ref{p7.2} that  $\int_{X^*\setminus (M^*)^n(U^*)}|\widehat f_k|^2\,d\mu=0$. Hence, for any  $\varepsilon>0$ and  for large enough $k$, using~\eqref{12.1}, we have
$$
\bigg(\int\limits_{X^*\setminus (M^*)^n(U^*)}|\widehat f|^2\,d\mu\bigg)^{1/2}\le \bigg(\int\limits_{X^*\setminus (M^*)^n(U^*)}|\widehat f-\widehat f_k|^2\,d\mu\bigg)^{1/2}+\bigg(\int\limits_{X^*\setminus (M^*)^n(U^*)}|\widehat f_k|^2\,d\mu\bigg)^{1/2}
$$
$$
\le \bigg(\int\limits_{X^*}|\widehat f-\widehat f_k|^2\,d\mu\bigg)^{1/2}=\bigg(\frac{\mu (U^*)}{\mu (U)}\int\limits_{X}| f- f_k|^2\,d\mu\bigg)^{1/2}
	<\varepsilon.
$$
Thus, $\supp \widehat f\subset  (M^*)^n(U^*)$, which completes the proof. 

{\it Sufficiently.} Let $\supp\widehat f \subset (M^*)^n(U^*)$. Then  $\widehat f\in L^1(X)$, and hence,  it follows from Proposition~\ref{p7.2} that the function  $F(f):={\widehat{(\overline{\widehat f})}}$ belongs to the space $S_n(X)$. 
Due to Proposition~\ref{p99.2} and Plancherel theorem, for every  $\varepsilon>0$  there exists $g\in \mathcal S $ such that 
$$
\|f-g\|_2<\varepsilon,\quad \|F(f)-F(g)\|_2<\varepsilon, 
$$
and  $F(g)=\overline g $ because of Proposition~\ref{p9.2}, that yields   $\|\overline f-F(f)\|_2=0$, i.e., the functions $\overline f$ и  $F(f)$ are equivalent. 
 \  $\Diamond$

\medskip

\begin{prop}\label{100}
	If $f\in L^2(X)$ and $\gamma\in H$, then 
\begin{equation}
   \widehat{f(\cdot\oplus \gamma)}(\omega)=\widehat f(\omega) \chi(\gamma,\omega), \quad \omega\in X^*.	
\label{100.1}
\end{equation}

\end{prop}

\noindent
{\bf Proof.} First we suppose that $f\in L(X)$. Since 
$$
    \chi (x,\omega) = \chi ((x\oplus \gamma)\ominus \gamma,\omega) = \chi (x\oplus \gamma,\omega)\overline{\chi (\gamma,\omega)},  
$$
for all $x\in X$ and $\omega\in X^*$, we have
$$
\widehat{f(\cdot\oplus \gamma)}(\omega)=\frac{1}{\mu(U)}\int\limits_{X}f(x\oplus\gamma)\overline{\chi (x,\omega)}d\mu(x)
$$
$$
=\frac{\chi(\gamma,\omega)}{\mu(U)}\int\limits_{X}f(x\oplus\gamma)\overline{\chi (x\oplus\gamma,\omega)}d\mu(x)
$$
$$
=\frac{\chi(\gamma,\omega)}{\mu(U)}\sum\limits_{h\in H}\int\limits_{U}f((x+h)\oplus\gamma)\overline{\chi ((x+h)\oplus \gamma,\omega)}d\mu(x).
$$
Taking into account that 
\begin{equation}
	\label{777}
		\{h'=h\oplus \gamma:\  h\in H\}=H\ \ {\rm and}\ \ 
	x+h=x\oplus h, \ \  (x+h)\oplus \gamma=x\oplus (h\oplus\gamma)
\end{equation}
for any $x\in U$ and $h\in H$, we obtain 
$$
\widehat{f(\cdot\oplus \gamma)}(\omega)
=\frac{\chi(\gamma,\omega)}{\mu(U)}\sum\limits_{h'\in H}\int\limits_{U}f(x\oplus h')\overline{\chi (x\oplus h', \omega)}d\mu(x). 
$$
Since here $x\oplus h'=x+h'$,  equality~\eqref{100.1} is proved for $f\in L(X)$.

Let now $f\in L^2(X)$, and let $f_n$  be the restriction of $f$ onto the set $M^n(U)$. Given $\varepsilon>0$,  we choose  $N$ such that
$\|f-f_N\|_2^2<\varepsilon$. Using~\eqref{777}, we get 
$$
\int\limits_{X}|f(x\oplus\gamma)-f_N(x\oplus\gamma)|^2d\mu(x)
=\sum\limits_{h\in H}\int\limits_{U}|f((x+h)\oplus\gamma)-f_N((x+h)\oplus\gamma)|^2d\mu(x)
$$
$$
=\sum\limits_{h'\in H}\int\limits_{U+h'}|f(x)-f_N(x)|^2d\mu(x)
=\int\limits_{X}|f(x)-f_N(x)|^2d\mu(x)<\varepsilon.
$$
Since $f_N\in L(X)$,  and the statement is proved already for such a function,  using the Plancherel theorem, we conclude 
$$
\|\widehat{f(\cdot\oplus \gamma)}-\widehat f (\cdot)\chi(\cdot, \gamma)\|_2\le
\|\widehat f-\widehat f_N \|_2+
	\|\widehat{f(\cdot\oplus \gamma)}-\widehat{f_N(\cdot\oplus \gamma)}\|_2
$$ 
$$
	=\frac{\mu(U^*)}{\mu(U)}(\|f-f_N\|_2+\|f(\cdot\oplus\gamma)-f_N(\cdot\oplus\gamma)\|_2)<C\varepsilon,
$$ 
which implies the required equality. \ \  $\Diamond$

\medskip 

\subsection*{Funding}
The second author was supported by the Russian Science Foundation (grant~
23-11-00178); Section~4 and Theorem~11 are due to this author.

\subsection*{Competing interests}
The authors have no competing interests to declare.

\end{document}